\documentclass[12pt,a4paper]{article}
\usepackage{amsfonts}
\usepackage{amsmath}
\usepackage{amsthm}
\usepackage[arrow,matrix]{xy}
\let\mod=\undefined
\DeclareMathOperator{\End}{End}
\DeclareMathOperator{\Ext}{Ext}
\DeclareMathOperator{\GL}{GL}
\DeclareMathOperator{\Hom}{Hom}
\DeclareMathOperator{\im}{Im}
\DeclareMathOperator{\ind}{ind}
\DeclareMathOperator{\Ker}{Ker}
\DeclareMathOperator{\mod}{mod}
\DeclareMathOperator{\rad}{rad}
\DeclareMathOperator{\Reg}{Reg}
\DeclareMathOperator{\rep}{rep}
\DeclareMathOperator{\Sing}{Sing}
\newcommand{\BA}{{\mathbb A}}
\newcommand{\BB}{{\mathbb B}}
\newcommand{\BM}{{\mathbb M}}
\newcommand{\BN}{{\mathbb N}}
\newcommand{\BZ}{{\mathbb Z}}
\newcommand{\CE}{{\mathcal E}}
\newcommand{\CF}{{\mathcal F}}
\newcommand{\CO}{{\mathcal O}}
\newcommand{\CT}{{\mathcal T}}
\newcommand{\CX}{{\mathcal X}}
\newcommand{\CY}{{\mathcal Y}}
\newcommand{\CZ}{{\mathcal Z}}
\newcommand{\dd}{{\mathbf d}}
\newcommand{\ov}{\overline}
\newcommand{\un}{\underline}
\newcommand{\bsmatrix}[1]{\left[\begin{smallmatrix} #1%
 \end{smallmatrix}\right]}
\newcommand{\psmatrix}[1]{\left(\begin{smallmatrix} #1%
 \end{smallmatrix}\right)}
\newtheorem{thm}{Theorem}[section]
\newtheorem{cor}[thm]{Corollary}
\newtheorem{lem}[thm]{Lemma}
\newtheorem{prop}[thm]{Proposition}
\newtheorem{step}{Step}
\numberwithin{equation}{section}
\begin{document}
\title{Orbit closures for representations of Dynkin quivers are
regular in codimension two
 }
\author{Grzegorz Zwara}
\date{\today}
\maketitle

\begin{abstract}
We develop reductions for classifications of singularities
of orbit closures in module varieties.
Then we show that the orbit closures for representations of
Dynkin quivers are regular in codimension two.
\end{abstract}

\section{Introduction and the main results}

Throughout the paper, $k$ denotes an algebraically closed field,
$A$ denotes a finitely generated associative $k$-algebra with
identity, and by a module we mean a left $A$-module whose
underlying $k$-space is finite dimensional.
Let $d$ be a positive integer and denote by $\BM_d(k)$
the algebra of $d\times d$-matrices with coefficients in $k$.
For an algebra $A$, the set $\mod_A(d)$ of algebra homomorphisms
$A\to\BM_d(k)$ has a natural structure of an affine variety.
Indeed, if $A\simeq k\langle X_1,\ldots,X_t\rangle/I$ for some
two-sided ideal $I$, then $\mod_A(d)$ can be identified with
the closed subset of $(\BM_d(k))^t$ given by vanishing of the
entries of all matrices $\rho(X_1,\ldots,X_t)$, $\rho\in I$.
Moreover, the general linear group $\GL(d)$ acts on $\mod_A(d)$
by conjugations
$$
g\star(M_1,\ldots,M_t)=(gM_1g^{-1},\ldots,gM_tg^{-1}),
$$
and the $\GL(d)$-orbits in $\mod_A(d)$ correspond bijectively
to the isomorphism classes of $d$-dimensional modules.
We shall denote by $\CO_M$ the $\GL(d)$-orbit in $\mod_A(d)$
corresponding to a $d$-dimensional module $M$.
An interesting problem is to study geometric properties of the
Zariski closure $\ov{\CO}_M$ of an orbit $\CO_M$ in $\mod_A(d)$.
We refer to \cite{BB}, \cite{BZ2}, \cite{Bgeo}, \cite{Bmin},
\cite{Bext}, \cite{SZder}, \cite{Zsmo}, \cite{Zuni} and
\cite{Zcodim1} for some results in this direction.

Following Hesselink (see~\cite[(1.7)]{Hes})
we call two pointed varieties $(\CX,x_0)$ and $(\CY,y_0)$
smoothly equivalent if there are smooth morphisms
$f:\CZ\to\CX$, $g:\CZ\to\CY$ and a point $z_0\in\CZ$
with $f(z_0)=x_0$ and $g(z_0)=y_0$.
This is an equivalence relation and the equivalence classes
will be denoted by $\Sing(\CX,x_0)$ and called the types
of singularities.
If $\Sing(\CX,x_0)=\Sing(\CY,y_0)$ then the variety $\CX$
is regular (respectively, normal, Cohen-Macaulay) at $x_0$
if and only if the same is true for the variety $\CY$ at $y_0$
(see~\cite[Section~17]{EGA} for more information about smooth
morphisms).
Obviously the regular points of the varieties give one type
of singularity, which we denote by $\Reg$.
Let $M$ and $N$ be $d$-dimensional modules with
$\CO_N\subseteq\ov{\CO}_M$, i.e., $N$ is a degeneration of $M$.
We shall write $\Sing(M,N)$ for $\Sing(\ov{\CO}_M,n)$, where $n$
is an arbitrary point of $\CO_N$.
It was shown recently (\cite[Theorem 1.1]{Zcodim1}) that
$\Sing(M,N)=\Reg$ provided $\dim\CO_M-\dim\CO_N=1$.
In this paper we investigate $\Sing(M,N)$ when
$\dim\CO_M-\dim\CO_N=2$.
First we prove some auxiliary result.

\begin{thm} \label{aux1}
Let $M'$, $N'$ and $X$ be modules such that
$\CO_{N'\oplus X}\subset\ov{\CO}_{M'\oplus X}$ and
$\dim\CO_{M'\oplus X}-\dim\CO_{N'\oplus X}=2$.
Then $\CO_{N'}\subset\ov{\CO}_{M'}$ and one of the following
cases holds:
\begin{enumerate}
\item[(1)] $\dim\CO_{M'}-\dim\CO_{N'}=1$ and
 $\Sing(M'\oplus X,N'\oplus X)=\Reg$;
\item[(2)] $\dim\CO_{M'}-\dim\CO_{N'}=2$ and
 $\Sing(M'\oplus X,N'\oplus X)=\Sing(M',N')$.
\end{enumerate}
\end{thm}

This allows to restrict our attention only to the case when
the modules $M$ and $N$ have no nonzero direct summands in common.
We shall say that such modules are disjoint.
We denote by $s(L)$ the number of summands in a decomposition
of a module $L$ into a direct sum of indecomposable modules.
The next result give us a further reduction for the problem of
description of the type $\Sing(M,N)$.

\begin{thm} \label{aux2}
Let $M$ and $N$ be disjoint modules such that
$\CO_N\subset\ov{\CO}_M$ and $\dim\CO_M-\dim\CO_N=2$.
Then $\Sing(M,N)=Reg$ if $s(N)\geq 3$.
\end{thm}

If $A=k[\varepsilon]/(\varepsilon^2)$, $M={}_AA$ and $N$ is
a direct sum of two simple $A$-modules, then $s(N)=2$,
$\CO_N\subset\ov{\CO}_M$, $\dim\CO_M-\dim\CO_N=2$ and
$$
\Sing(M,N)=\Sing\left(\left\{\bsmatrix{x&y\\ z&-x};\;x^2+yz=0
\right\},\bsmatrix{0&0\\ 0&0}\right)
$$
is the type of Kleinian singularity $\BA_2$.
Hence orbit closures in module varieties may be singular in
codimension two even for very simple algebras.
However this is not true for the modules over the path algebras
of Dynkin quivers.
We add that Theorems~\ref{aux1} and \ref{aux2} are used in
the proof of our main result stated below.

\begin{thm} \label{main}
Let $M$ be a module over the path algebra of a Dynkin quiver.
Then the variety $\ov{\CO}_M$ is regular in codimension two.
\end{thm}

Let $Q=(Q_0,Q_1,s,e)$ be a finite quiver.
Here $Q_0$ is a finite set of vertices, $Q_1$ is a finite set
of arrows, and $s,e:Q_1\to Q_0$ are functions such that any arrow
$\alpha\in Q_1$ has the starting vertex $s(\alpha)$ and the
ending vertex $e(\alpha)$.
Let $\dd=(d_i)_{i\in Q_0}\in\BN^{Q_0}$.
We define the vector space
$$
\rep_Q (\dd)=\prod_{\alpha\in Q_1}\BM_{d_{e(\alpha)}\times
d_{s(\alpha)}}(k),
$$
where $\BM_{d'\times d''}(k)$ denotes the set of
$d'\times d''$-matrices with coefficients in $k$ for any
$d',d''\in\BN$.
The product $\GL(\dd)=\prod_{i\in Q_0}\GL(d_i)$ of general linear
groups acts on $\rep_Q(\dd)$ via
$$
g\star V=(g_{e(\alpha)}V_\alpha g_{s(\alpha)}^{-1})
_{\alpha\in Q_1},
$$
for any $g=(g_i)_{i\in Q_0}\in\GL(\dd)$ and
$V=(V_\alpha)_{\alpha\in Q_1}\in\rep_Q(\dd)$.
Using an equivalence described by Bongartz in \cite{Bgeo} we
can reformulate Theorem~\ref{main} as follows.

\begin{cor} \label{cormain}
Let $Q$ be a Dynkin quiver and $\dd\in\BN^{Q_0}$.
Then the closures of the $\GL(\dd)$-orbits in $\rep_Q(\dd)$ are
regular in codimension two.
\end{cor}

Let $Q:1\xrightarrow{\alpha}2$ be a Dynkin quiver of type $\BA_2$
and $\dd=(2,2)\in\BN^{Q_0}$.
Then $\rep_Q(\dd)=\BM_{2\times 2}(k)$ and the orbit closure
$$
\ov{\GL(\dd)\star\bsmatrix{1&0\\ 0&0}}
=\left\{\bsmatrix{x&y\\ z&t};\;xt-yz=0\right\}
$$
is a singular variety of dimension three.
This shows that ``codimension two'' in Corollary~\ref{cormain}
(and in Theorem~\ref{main}) cannot be improved by
``codimension three''.

We shall consider in Section~\ref{degen} some properties
of short exact sequences, dimensions of homomorphism spaces
and degenerations of modules.
Section~\ref{smooth} contains some sufficient conditions on
regularity of $\Sing(M,N)$.
Sections~\ref{proofaux1}, \ref{proofaux2} and \ref{proofmain}
are devoted to the proofs of Theorems~\ref{aux1}, \ref{aux2}
and \ref{main}, respectively.

For basic background on the representation theory of algebras and
quivers we refer to \cite{ARS} and \cite{Rin}.
The author gratefully acknowledges support from the Polish
Scientific Grant KBN No.\ 1 P03A 018 27.

\section{Degenerations of modules}
\label{degen}

Let $\mod A$ denote the category of finite dimensional left
$A$-modules and $\rad(\mod A)$ denote the Jacobson radical of
the category $\mod A$.
We can describe $\rad(\mod A)$ as the two-sided ideal of $\mod A$
generated by nonisomorphisms between indecomposable modules.
We abbreviate by $[X,Y]$ the dimension $\dim_k\Hom_A(X,Y)$ for
any modules $X$ and $Y$.
Recall that by a module we mean an object of $\mod A$.

\begin{lem} \label{formulacodim}
Let $M$ and $N$ be modules with $\dim_kM=\dim_kN$.
Then $\dim\CO_M-\dim\CO_N=[N,N]-[M,M]$.
\end{lem}

\begin{proof}
Let $L$ be a $d$-dimensional module and choose a point $l$ in
$\CO_L$.
Since the isotropy group of $l$ can be identified with the
group of $A$-automorphisms of $L$ and the latter is a nonempty
and open subset of the vector space $\End_A(L)$, then
we conclude the formula
$$
\dim\CO_L=\dim\GL(d)-[L,L].
$$
We get the claim by applying the formula for $L=M$ and $L=N$.
\end{proof}

We shall need the following three simple facts on short exact
sequences.

\begin{lem} \label{seqineq}
Let $X$ be a module and
$\sigma:0\to U\xrightarrow{f}W\xrightarrow{g}V\to 0$ be an exact
sequence in $\mod A$.
Then:
\begin{enumerate}
\item[(1)] $\delta_\sigma(X):=[U\oplus V,X]-[W,X]\geq 0$ and the
equality holds if and  only if any homomorphism in $\Hom_A(U,X)$
factors through $f$;
\item[(2)] $\delta'_\sigma(X):=[X,U\oplus V]-[X,W]\geq 0$ and the
equality holds if and  only if any homomorphism in $\Hom_A(X,V)$
factors through $g$.
\end{enumerate}
\end{lem}

\begin{proof}
The claim follow from the induced exact sequences
\begin{align*}
&0\to\Hom_A(V,X)\xrightarrow{\Hom_A(g,X)}\Hom_A(W,X)
\xrightarrow{\Hom_A(f,X)}\Hom_A(U,X),\\
&0\to\Hom_A(X,U)\xrightarrow{\Hom_A(X,f)}\Hom_A(X,W)
\xrightarrow{\Hom_A(X,g)}\Hom_A(X,V).
\end{align*}
\end{proof}

\begin{lem} \label{whensplits}
Let $\sigma:\;0\to U\xrightarrow{f}W\xrightarrow{g}V\to 0$
be an exact sequence in $\mod A$.
Then the following conditions are equivalent.
\begin{enumerate}
\item[(1)] The sequence $\sigma$ splits.
\item[(2)] $W\simeq U\oplus V$.
\item[(3)] $\delta_\sigma(U)=0$.
\item[(4)] $\delta'_\sigma(V)=0$.
\end{enumerate}
\end{lem}

\begin{proof}
Clearly the condition \textit{(1)} implies \textit{(2)}, and the
condition \textit{(2)} implies \textit{(3)} and \textit{(4)}.
Applying Lemma~\ref{seqineq} we get that \textit{(3)} implies that
the endomorphism $1_U$ factors through $f$, which means that $f$
is a section and \textit{(1)} holds.
Similarly, it follows from \textit{(4)} that $g$ is a retraction
and \textit{(1)} holds.
\end{proof}

\begin{lem} \label{howtosplitiso}
Let
$$
0\to U\xrightarrow{\psmatrix{f_1\\ f_2}}W_1\oplus W_2
\xrightarrow{\psmatrix{g_{1,1}&g_{1,2}\\ g_{2,1}&g_{2,2}}}
V_1\oplus V_2\to 0
$$
be an exact sequence in $\mod A$ such that $g_{1,1}$ is an
isomorphism.
Then
$$
0\to U\xrightarrow{f_2}W_2\xrightarrow{g'}V_2\to 0.
$$
is also an exact sequence in $\mod A$, where
$g'=g_{2,2}-g_{2,1}g_{1,1}^{-1}g_{1,2}$.
\end{lem}

\begin{proof}
Straightforward.
\end{proof}

The next result follows from \cite[Theorem 1.1]{Zgiv} and from
Lemma~\ref{howtosplitiso} and its dual.

\begin{thm} \label{specialseq}
Let $M$ and $N$ be modules.
Then the inclusion $\CO_N\subseteq\ov{\CO}_M$ is equivalent
to each of the following conditions:
\begin{enumerate}
\item[(1)] There is an exact sequence
$0\to Z\xrightarrow{f}Z\oplus M\xrightarrow{g}N\to 0$ in $\mod A$
for some module $Z$.
\item[(2)] There is an exact sequence
$0\to N\xrightarrow{f'}M\oplus Z'\xrightarrow{g'}Z'\to 0$ in
$\mod A$ for some module $Z'$.
\end{enumerate}
Moreover, we may assume that $f$ and $g'$ belong to
$\rad(\mod A)$.
\end{thm}

\begin{cor} \label{UMVdeg}
Let
$$
\sigma:\quad 0\to U\to M\to V\to 0
$$
be an exact sequence in $\mod A$.
Then $\CO_{U\oplus V}\subseteq\ov{\CO}_M$.
\end{cor}

\begin{proof}
We apply Theorem~\ref{specialseq} to a direct sum of $\sigma$ and
the exact sequence $0\to 0\to U\xrightarrow{1_U}U\to 0$.
\end{proof}

\begin{lem} \label{deghom}
Let $M$ and $N$ be modules such that $\CO_N\subseteq\ov{\CO}_M$.
Then
$$
\delta_{M,N}(X):=[N,X]-[M,X]\geq 0\quad\text{and}\quad
\delta'_{M,N}(X):=[X,N]-[X,M]\geq 0
$$
for any module $X$.
\end{lem}

\begin{proof}
We get an exact sequence $\sigma: 0\to Z\to Z\oplus M\to N\to 0$
in $\mod A$, by Theorem~\ref{specialseq}. Then the claim follows
from Lemma~\ref{seqineq} and the equalities
$\delta_{M,N}(X)=\delta_\sigma(X)$ and
$\delta'_{M,N}(X)=\delta'_\sigma(X)$ for any module $X$.
\end{proof}

Let $M$ and $N$ be modules with $\CO_N\subseteq\ov{\CO}_M$ and
$\sigma$ be a short exact sequence in $\mod A$.
We shall use frequently without refereing the following
obvious properties of the nonnegative integers $\delta(L)$:
\begin{itemize}
\item $\delta(X)=\delta(Y)$ if $X\simeq Y$,
\item $\delta(X\oplus Y)=\delta(X)+\delta(Y)$,
\item $\delta(X\oplus Y)=0$ implies $\delta(X)=0$,
\end{itemize}
where $X$ and $Y$ are modules and $\delta$ is an abbreviation of
$\delta_\sigma$, $\delta'_\sigma$, $\delta_{M,N}$ or
$\delta'_{M,N}$.

\section{Smooth points of orbit closures}
\label{smooth}

Throughout the section let $M$ and $N$ be $d$-dimensional
modules such that $\CO_N\subseteq\ov{\CO}_M$, and let $\CF_{M,N}$
and $\CF'_{M,N}$ denote complete sets of pairwise nonisomorphic
modules $X$ such that $\delta_{M,N}(X)=0$ and
$\delta'_{M,N}(X)=0$, respectively.

Let $U,V\in\mod A$.
We denote by $\BZ^1_A(V,U)$ the group of cocycles, i.e., the
$k$-linear maps $Z:A\to\Hom_k(V,U)$ satisfying
$$
Z(aa')=Z(a)V(a')+U(a)Z(a'),\qquad\text{for all }a,a'\in A.
$$
The group $\BZ^1_A(V,U)$ contains the group of coboundaries
$$
\BB^1_A(V,U)=\{hV-Uh;\;h\in\Hom_k(V,U)\}.
$$
This leads to the $k$-functor
$$
\BZ^1_A(-,-):\mod A\times\mod A\to\mod k
$$
and its $k$-subfunctor $\BB^1_A(-,-)$.
Any cocycle $Z$ in $\BZ^1_A(V,U)$ induces an exact sequence
$$
\sigma_Z: 0\to U\xrightarrow{\alpha_Z}W_Z\xrightarrow{\beta_Z}
V\to 0
$$
in $\mod A$.
Then the cocycle $Z$ is a coboundary if and only if the sequence
$\sigma_Z$ splits, which is equivalent to the fact that
$W_Z\simeq U\oplus V$, by Lemma~\ref{whensplits}.
Let
\begin{align*}
\CZ_{M,N}(V,U)=\{Z\in\BZ^1_A(V,U);\quad
 &\delta_{\sigma_Z}(X)=0\text{ for any }X\in\CF_{M,N},\\
&\;\delta'_{\sigma_Z}(Y)=0\text{ for any }Y\in\CF'_{M,N}\}.
\end{align*}
Obviously $\CZ_{M,N}(V,U)$ contains $\BB^1_A(V,U)$ and does not
depend on the choice of representatives of isomorphism classes
of modules in the definition of the sets $\CF_{M,N}$ and
$\CF'_{M,N}$.

\begin{lem} \label{reformCZ}
A cocycle $Z\in\BZ^1_A(V,U)$ belongs to $\CZ_{M,N}(V,U)$ if and
only if
$$
\BZ^1_A(V,f)(Z)\in\BB^1_A(V,X)\quad\text{and}\quad
\BZ^1_A(g,U)(Z)\in\BB^1_A(Y,U)
$$
for any modules $X\in\CF_{M,N}$, $Y\in\CF'_{M,N}$ and any
homomorphisms $f:U\to X$, $g:Y\to V$.
\end{lem}

\begin{proof}
Let $Z$ be a cocycle in $\BZ^1_A(V,U)$.
By duality, it suffices to show that $\delta_{\sigma_Z}(X)=0$
if and only if the cocycle $\BZ^1_A(V,f)(Z)$ is a coboundary
for any homomorphism $f:U\to X$.
By Lemma~\ref{seqineq}, the equality $\delta_{\sigma_Z}(X)=0$
means that any homomorphism in $\Hom_A(U,X)$ factors through
$\alpha_Z$.
Let $Z'=\BZ^1_A(V,f)(Z)$ for some homomorphism $f:U\to X$.
We consider the pushout of $\sigma_Z$ under $f$:
$$
\xymatrix{
\sigma_Z:&0\ar[r]&U\ar[r]^{\alpha_Z}\ar[d]_f
 &W_Z\ar[r]^{\beta_Z}\ar[d]&V\ar[r]\ar@{=}[d]&0\\
\sigma_{Z'}:&0\ar[r]&X\ar[r]^{\alpha_{Z'}}
 &W_{Z'}\ar[r]^{\beta_{Z'}}&V\ar[r]&0.
}
$$
Then $f$ factors through $\alpha_Z$ if and only if the sequence
$\sigma_{Z'}$ splits, and the latter means that the cocycle $Z'$
is a coboundary.
\end{proof}

\begin{lem} \label{Zsub}
$\CZ_{M,N}(-,-)$ is a $k$-subfunctor of $\BZ^1_A(-,-)$.
\end{lem}

\begin{proof}
Let $U$ and $V$ be modules.
We take $X\in\CF_{M,N}$ and $Y\in\CF'_{M,N}$.
Then $\CZ_{M,N}(V,U)$ is a $k$-space, by Lemma~\ref{reformCZ}
and since the appropriate maps $\BZ^1_A(V,f)$ and $\BZ^1_A(g,U)$
are $k$-linear.
Let $Z$ be a cocycle in $\CZ_{M,N}(V,U)$.
We set $Z'=\BZ^1_A(V,f')(Z)$,
where $f':U\to U'$ is a homomorphism for some module $U'$.
Then
$$
\BZ^1_A(V,\tilde{f})(Z')=\BZ^1_A(V,\tilde{f}f')(Z)\in\BB^1_A(V,X)
$$
for any homomorphism $\tilde{f}:U'\to X$ and
\begin{align*}
\BZ^1_A(\tilde{g},U')(Z')&=\BZ^1_A(\tilde{g},f')(Z)
 =\BZ^1_A(Y,f')\left(\BZ^1_A(\tilde{g},U)(Z)\right)\\
&\in\BZ^1_A(Y,f')\left(\BB^1_A(Y,U)\right)\subseteq\BB^1_A(Y,U')
\end{align*}
for any homomorphism $\tilde{g}:Y\to V$.
This shows that the cocyle $Z'$ belongs to $\CZ_{M,N}(V,U')$.
Dually the cocycle $\BZ^1_A(g',U)(Z)$ belongs to $\CZ_{M,N}(V',U)$
for any module $V'$ and any homomorphism $g':V'\to V$.
\end{proof}

The module variety $\mod_A(d)$ is the underlying variety of
an affine $k$-scheme $\un{\mod}_A^d$ of finite type, which
represents the functor
$$
\mod_A^d:(Commutative\;k\text{--}algebras)\to(Sets),
$$
where $\mod_A^d(R)$ is the set of $k$-algebra homomorphisms
from $A$ to the algebra of $d\times d$-matrices with coefficients
in a commutative $k$-algebra $R$ (see \cite{Bgeo}, \cite{Gopen}).
We denote by $\CT_{\CX,x}$ the tangent space of a $k$-scheme
$\CX$ at a point $x$.
Let $n$ be a (closed) point of $\CO_N$.
Then the tangent space $\CT_{\un{\mod}_A^d,n}$ corresponds to
the preimage of $n$ via the canonical map
$$
\mod_A^d(k[\varepsilon]/(\varepsilon^2))\to
\mod_A^d(k),
$$
and the latter corresponds to the group of cocycles
$\BZ^1_A(N,N)$.
Hence we get a canonical $k$-isomorphism
$$
\Phi:\CT_{\un{\mod}_A^d,n}\xrightarrow{\simeq}\BZ^1_A(N,N).
$$
Furthermore, $\Phi(\CT_{\CO_N,n})=\BB^1_A(N,N)$ which gives
the isomorphism
$$
\ov{\Phi}:\CT_{\un{\mod}_A^d,n}/\CT_{\CO_N,n}
\xrightarrow{\simeq}\Ext^1_A(N,N)
$$
known as a Voigt result (see \cite[Proposition 1.1]{Gopen}).
Here and later on, the group $\Ext_A^1(V,U)$ of extensions
of $V$ by $U$ is identified with the quotient
$\BZ^1_A(V,U)/\BB^1_A(V,U)$ for any modules $U$ and $V$.

\begin{lem} \label{tanginCZ}
Let $n\in\CO_N$. Then
$\Phi\left(\CT_{\ov{\CO}_M,n}\right)\subseteq\CZ_{M,N}(N,N)$.
\end{lem}

\begin{proof}
We have to recall some notation and results of Section 3 in
\cite{Zsmo} (see also the proof of \cite[Proposition 2.2]{Zuni}).
Let $X$ be a module and
$$
\mod_{A,X,t}^d:(Commutative\;k\text{--}algebras)\to(Sets)
$$
be the subfunctor of $\mod_A^d$ defined in \cite[(3.3)]{Zuni},
where $t=[X,M]$.
This functor is represented by an affine $k$-subscheme
$\CX=\un{\mod}_{A,X,t}^d$ of $\un{\mod}_A^d$ such that the
underlying variety is given by
$$
\CX_{red}=\{l\in\mod_A(d);\;[X,L]=t\}.
$$
Here $L$ denotes a module corresponding to a point $l$ in
$\mod_A(d)$.
Assume that $\delta'_{M,N}(X)=0$.
Then the orbits $\CO_M$ and $\CO_N$ are included in $\CX_{red}$.
Therefore $\CT_{\ov{\CO}_M,n}$ is contained in $\CT_{\CX,n}$.
On the other hand, the tangent space $\CT_{\CX,n}$ corresponds
to the preimage of $n$ via the canonical map
$$
\mod_{A,X,t}^d(k[\varepsilon]/(\varepsilon^2))\to
\mod_{A,X,t}^d(k).
$$
Furthermore, by \cite[Lemma 3.11]{Zsmo}, the latter corresponds
to the subset of $\BZ^1_A(N,N)$ consisting of the cocycles $Z$
such that $\delta'_{\sigma_Z}(X)=0$.
Hence $\Phi(\CT_{\ov{\CO}_M,n})$ is contained in
$$
\left\{Z\in\BZ^1_A(N,N);\quad\delta'_{\sigma_Z}(X)=0
\text{ for any }X\in\CF'_{M,N}\right\}.
$$
By duality, $\Phi(\CT_{\ov{\CO}_M,n})$ is also contained in
$$
\left\{Z\in\BZ^1_A(N,N);\quad\delta_{\sigma_Z}(X)=0
\text{ for any }X\in\CF_{M,N}\right\},
$$
and the claim follows from the definition of $\CZ_{M,N}(N,N)$.
\end{proof}

We define the quotient
$\CE_{M,N}(V,U)=\CZ_{M,N}(V,U)/\BB^1_A(V,U)$ for any modules
$U$ and $V$.
An immediate consequence of Lemmas~\ref{reformCZ} and \ref{Zsub}
is the following fact.

\begin{cor} \label{Efunctor}
$\CE_{M,N}(-,-)$ is a $k$-subfunctor of
$$
\Ext^1_A(-,-):\mod A\times\mod A\to\mod k
$$
and
$$
\CE_{M,N}(V,U)
=\!\!\!\bigcap_{\substack{X\in\CF_{M,N}\\ f\in\Hom_A(U,X)}}
\!\!\!\Ker\left(\Ext^1_A(V,f)\right)\;\cap
\!\!\!\bigcap_{\substack{Y\in\CF'_{M,N}\\ g\in\Hom_A(Y,V)}}
\!\!\!\Ker\left(\Ext^1_A(g,U)\right)
$$
for any modules $U$ and $V$.
\end{cor}

Now we are ready to formulate our first sufficient conditions
for regularity of points in $\ov{\CO}_M$.

\begin{prop} \label{gencriterion}
$\dim_k\CE_{M,N}(N,N)\geq[N,N]-[M,M]$ and the equality implies
that $\Sing(M,N)=\Reg$.
\end{prop}

\begin{proof}
Let $n\in\CO_N$.
Combining Lemmas~\ref{formulacodim} and \ref{tanginCZ} we get
\begin{align*}
\dim_k\CE_{M,N}(N,N)&=\dim_k\CZ_{M,N}(N,N)-\dim_k\BB^1_A(N,N)\\
&\geq\dim_k\CT_{\ov{\CO}_M,n}-\dim_k\CT_{\CO_N,n}
 =\dim_k\CT_{\ov{\CO}_M,n}-\dim\CO_N\\
&\geq\dim\ov{\CO}_M-\dim\CO_N=[N,N]-[M,M].
\end{align*}
Moreover, the equality $\dim_k\CE_{M,N}(N,N)=[N,N]-[M,M]$ implies
that
$$
\dim_k\CT_{\ov{\CO}_M,n}=\dim\ov{\CO}_M,
$$
which means that $\Sing(M,N)=\Reg$, as the variety $\ov{\CO}_M$
is irreducible.
\end{proof}

As a consequence of the above proposition one can conclude the
following useful result (see \cite[Proposition 2.2]{Zuni}).

\begin{prop} \label{standcriterion}
Assume that one of the following cases holds.
\begin{enumerate}
\item[(1)] There is an exact sequence
$\sigma: 0\to Z\to Z\oplus M\to N\to 0$ in $\mod A$ and
$\delta'_{M,N}(Z\oplus M)=0$ for some module $Z$.
\item[(2)] There is an exact sequence
$\sigma': 0\to N\to M\oplus Z'\to Z'\to 0$ in $\mod A$ and
$\delta_{M,N}(M\oplus Z')=0$ for some module $Z'$.
\end{enumerate}
Then $\Sing(M,N)=\Reg$.
\end{prop}

\begin{proof}
\textit{(1)}.
We may assume that $Z\oplus M$ belongs to $\CF'_{M,N}$.
By Corollary~\ref{Efunctor}, $\CE_{M,N}(N,N)$ is contained in
the kernel of the last map in the following long exact sequence
induced by $\sigma$:
\begin{multline*}
0\to\Hom_A(N,N)\to\Hom_A(Z\oplus M,N)\to\Hom_A(Z,N)\to\\
\to\Ext^1_A(N,N)\to\Ext^1_A(Z\oplus M,N).
\end{multline*}
Consequently,
$$
\dim_k\CE_{M,N}(N,N)\leq\delta_\sigma(N)
=\delta_{M,N}(N)+\delta'_{M,N}(M)=[N,N]-[M,M].
$$
Hence the claim follows from Proposition~\ref{gencriterion}.

We proceed dually in case \textit{(2)}.
\end{proof}

\begin{cor} \label{corcriterion}
Let $\sigma:\;0\to U\to M\to V\to 0$ be an exact sequence in
$\mod A$ such that $\delta'_\sigma(U\oplus M)=0$ or
$\delta_\sigma(M\oplus V)=0$.
Then
$$
\Sing(M,U\oplus V)=\Reg.
$$
\end{cor}

\begin{proof}
If $\delta'_\sigma(U\oplus M)=0$ then it suffices to apply
Proposition~\ref{standcriterion} for $Z=U$ and the direct sum
of $\sigma$ and the sequence $0\to 0\to U\xrightarrow{1_U}U\to 0$.
We proceed in a similar way if $\delta_\sigma(M\oplus V)=0$.
\end{proof}

We conclude from the proof of \cite[Theorem 1.1]{Zcodim1} and
its dual the following result.

\begin{thm} \label{codim1}
Assume that $\dim\CO_M-\dim\CO_N=1$.
Then:
\begin{enumerate}
\item[(1)] $\delta_{M,N}(M)=\delta'_{M,N}(M)=0$ and
$\delta_{M,N}(N)=\delta'_{M,N}(N)=1$;
\item[(2)] there is an exact sequence
$0\to Z\to Z\oplus M\to N\to 0$ in $\mod A$ for some
indecomposable module $Z$ with $\delta'_{M,N}(Z)=0$;
\item[(3)] there is an exact sequence
$0\to N\to M\oplus Z'\to Z'\to 0$ in $\mod A$ for some
indecomposable module $Z'$ with $\delta_{M,N}(Z')=0$.
\end{enumerate}
In particular $\Sing(M,N)=\Reg$.
\end{thm}

\section{Reduction to disjoint modules}
\label{proofaux1}

Combining Lemmas~\ref{seqineq} and \ref{howtosplitiso} we get
the following fact.

\begin{lem} \label{splitdirectsum}
Let
$$
\sigma:\quad 0\to U\xrightarrow{f}W
\xrightarrow{g}V_1\oplus V_2\to 0
$$
be an exact sequence in $\mod A$ such that
$\delta'_\sigma(V_1)=0$.
Then $W=W_1\oplus W_2$ for some modules $W_1\simeq V_1$ and $W_2$
such that there is an exact sequence
$$
\eta:\quad 0\to U\xrightarrow{f'}W_2\xrightarrow{g'}V_2\to 0
$$
in $\mod A$ with $f':U\to W_2$ being a component of
$f:U\to W_1\oplus W_2$.
\end{lem}

We denote by $\mu(L,Y)$ the multiplicity of an indecomposable
module $Y$ as a direct summand of a module $L$.

\begin{lem} \label{ostra}
Let $M$ and $N$ be modules such that $\CO_N\subseteq\ov{\CO}_M$.
Let $Y$ be an indecomposable module such that $\mu(M,Y)<\mu(N,Y)$.
Then $\delta_{M,N}(Y)>0$ or $\delta'_{M,N}(Y)>0$.
\end{lem}

\begin{proof}
Applying Theorem~\ref{specialseq} we get an exact sequence
$$
\sigma:\quad 0\to Z\xrightarrow{f}Z\oplus M\to N\to 0
$$
in $\mod A$ such that $f$ belongs to $\rad(\mod A)$.
Let $Y$ be an indecomposable $A$-module such that
$p:=\mu(N,Y)>\mu(M,Y)$.
Assume that $\delta'_{M,N}(Y)=0$.
Then $\delta'_\sigma(Y^p)=\delta'_{M,N}(Y^p)=0$ and $Y^p$ is
isomorphic to a direct summand of $Z\oplus M$, by
Lemma~\ref{splitdirectsum}.
Therefore $\mu(Z\oplus M,Y)\geq p$ and consequently $\mu(Z,Y)>0$.
This means that there is a retraction $h:Z\to Y$.
We know that $h$ does not factor through $f$, as the latter
belongs to $\rad(\mod A)$.
Hence $\delta_{M,N}(Y)=\delta_\sigma(Y)>0$, by
Lemma~\ref{seqineq}.
\end{proof}

\begin{lem} \label{nnmm}
Let $M'$, $N'$ and $X$ be modules such that
$\CO_{N'\oplus X}\subset\ov{\CO}_{M'\oplus X}$ and
$M'\not\simeq N'$.
Then $[N',N']>[M',M']$.
\end{lem}

\begin{proof}
Let $M=M'\oplus X$ and $N=N'\oplus X$.
Since $M'$ and $N'$ are not isomorphic and $\dim_kM'=\dim_kN'$,
then there is an indecomposable $A$-module $Y$ such that
$\mu(N',Y)>\mu(M',Y)$, or equivalently, $\mu(N,Y)>\mu(M,Y)$.
Consequently $\delta_{M,N}(Y)>0$ or $\delta'_{M,N}(Y)>0$, by
Lemma~\ref{ostra}.
Therefore the claim follows from the inequalities
\begin{align*}
[N',N']-[M',M']&=\delta_{M,N}(N')+\delta'_{M,N}(M')\geq
 \delta_{M,N}(N')\geq\delta_{M,N}(Y),\\
[N',N']-[M',M']&=\delta'_{M,N}(N')+\delta_{M,N}(M')\geq
 \delta'_{M,N}(N')\geq\delta'_{M,N}(Y).
\end{align*}
\end{proof}

We shall need the following cancellation properties proved
by Bongartz (see \cite[Corollary 2.5]{Bext} and
\cite[Theorem 2]{Bmin}).

\begin{thm} \label{cancel}
Let $M'$, $N'$ and $X$ be modules such that
$\CO_N\subseteq\ov{\CO}_M$ for $M=M'\oplus X$ and $N=N'\oplus X$.
\begin{enumerate}
\item[(1)] If $\delta_{M,N}(X)=0$ or $\delta'_{M,N}(X)=0$ then
$\CO_{N'}\subseteq\ov{\CO}_{M'}$.
\item[(2)] If $\delta_{M,N}(X)=0$ and $\delta'_{M,N}(X)=0$ then
$\Sing(M,N)=\Sing(M',N')$.
\end{enumerate}
\end{thm}
\bigskip

\noindent
\textit{Proof of Theorem~\ref{aux1}.}
Let $M'$, $N'$ and $X$ be modules such that
$\CO_N\subset\ov{\CO}_M$ and $\dim\CO_M-\dim\CO_N=2$,
where $M=M'\oplus X$ and $N=N'\oplus X$.
In particular, the modules $M'$ and $N'$ are not isomorphic and
$$
2=[N,N]-[M,M]=([N',N']-[M',M'])+\delta_{M,N}(X)+\delta'_{M,N}(X).
$$
On the other hand $[N',N']-[M',M']\geq 1$, by Lemma~\ref{nnmm}.
Therefore
$$
\dim\CO_{M'}-\dim\CO_{N'}=[N',N']-[M',M']\in\{1,2\},
$$
and at least one of the numbers $\delta_{M,N}(X)$ and
$\delta'_{M,N}(X)$ is zero.
Consequently $\CO_{N'}\subseteq\ov{\CO}_{M'}$, by
Theorem~\ref{cancel}.

We first consider the case $\dim\CO_{M'}-\dim\CO_{N'}=1$.
By duality, we may assume that $\delta'_{M,N}(X)=0$.
Using Theorem~\ref{codim1} we derive the exact sequence
$$
\sigma:\; 0\to Z\to Z\oplus M'\to N'\to 0
$$
in $\mod A$ for some module $Z$ such that
$\delta'_{M',N'}(Z\oplus M')=0$.
Hence
$$
\delta'_{M,N}(Z\oplus M)=\delta'_{M,N}(Z\oplus M')
+\delta'_{M,N}(X)=\delta'_{M',N'}(Z\oplus M')=0.
$$
Let
$$
0\to Z\to Z\oplus M\to N\to 0
$$
be a direct sum of $\sigma$ and the short exact sequence
$$
0\to 0\to X\xrightarrow{1_X}X\to 0.
$$
Then $\Sing(M,N)=\Reg$, by Proposition~\ref{standcriterion}.

It remains to consider the case $\dim\CO_{M'}-\dim\CO_{N'}=2$.
Then $\delta_{M,N}(X)=\delta'_{M,N}(X)=0$.
Hence $\Sing(M,N)=\Sing(M',N')$, by Theorem~\ref{cancel}.
\qed

\section{Reduction to at most two summands}
\label{proofaux2}

We shall need the following result which can be derived from
the proof of \cite[Theorem 2.3]{Zlik}.

\begin{prop} \label{tubes}
Let $0\to Z\xrightarrow{f}Z\oplus M\to N\to 0$ be an exact
sequence in $\mod A$ such that the homomorphism $f$ belongs to
$\rad(\mod A)$.
Then there are a positive integer $h$ and exact sequences
$$
\sigma_i:\quad 0\to N_i\to N_{i-1}\oplus N_{i+1}\to N_i\to 0,
\qquad i=1,2,\ldots,h,
$$
in $\mod A$ for some modules $N_0,N_1,\ldots,N_{h+1}$ such that
$N_0=0$, $N_1\simeq N$, $N_{h+1}\simeq N_h\oplus M$ and
$Z$ is isomorphic to a direct summand of $N_h$.
\end{prop}

\begin{lem} \label{degmix}
Let $0\to Z\xrightarrow{f}Z\oplus M\to N\to 0$ be an exact
sequence in $\mod A$ such that $f$ belongs to $\rad(\mod A)$.
Let $\tilde{M}$ and $\tilde{N}$ be modules such that
$\CO_{\tilde{N}}\subseteq\ov{\CO}_{\tilde{M}}$ and
$\delta_{M,N}(\tilde{M})=\delta_{M,N}(\tilde{N})
=\delta'_{\tilde{M},\tilde{N}}(N)=0$.
Then $\delta'_{\tilde{M},\tilde{N}}(Z)=0$.
\end{lem}

\begin{proof}
We use Proposition~\ref{tubes} and the notation introduced
there.
Then
$$
\sum_{i=1}^h \delta_{\sigma_i}(\tilde{M})=\delta_{M,N}(\tilde{M})
=0.
$$
This implies that
$$
2\cdot[N_i,\tilde{M}]-[N_{i+1},\tilde{M}]-[N_{i-1},\tilde{M}]
=\delta_{\sigma_i}(\tilde{M})=0,\qquad i=1,2,\ldots,h.
$$
Proceeding by induction on $i$, one can show that
$$
[N_i,\tilde{M}]=i\cdot[N,\tilde{M}],\qquad i=0,1,\ldots,h+1.
$$
In a similar way we get
$$
[N_i,\tilde{N}]=i\cdot[N,\tilde{N}],\qquad i=0,1,\ldots,h+1.
$$
In particular
$$
\delta'_{\tilde{M},\tilde{N}}(N_h)
=h\cdot\delta'_{\tilde{M},\tilde{N}}(N)=0\quad\text{and}\quad
\delta'_{\tilde{M},\tilde{N}}(Z)=0,
$$
as $Z$ is isomorphic to a direct summand of $N_h$.
\end{proof}

\begin{prop} \label{twodeg}
Let $M'$, $M''$, $N'$ and $N''$ be modules such that
$M'\not\simeq N'$, $M''\not\simeq N''$,
$\CO_{N'}\subset\ov{\CO}_{M'}$,
$\CO_{N''}\subset\ov{\CO}_{M''}$
and
$$
\dim\CO_{M'\oplus M''}-\dim\CO_{N'\oplus N''}=2.
$$
Then $\Sing(M'\oplus M'',N'\oplus N'')=\Reg$.
\end{prop}

\begin{proof}
It follows from the assumptions and Lemma~\ref{formulacodim}
that the integers
$$
c'=\dim\CO_{M'}-\dim\CO_{N'}\quad\text{and}\quad
c''=\dim\CO_{M''}-\dim\CO_{N''}
$$
are positive and
\begin{align*}
2&=c'+c''+\delta_{M',N'}(N'')+\delta'_{M',N'}(N'')
 +\delta_{M'',N''}(M')+\delta'_{M'',N''}(M')\\
&=c'+c''+\delta_{M',N'}(M'')+\delta'_{M',N'}(M'')
 +\delta_{M'',N''}(N')+\delta'_{M'',N''}(N').
\end{align*}
Hence $c'=c''=1$ and
\begin{equation} \label{1outof3}
\begin{aligned}
&\delta_{M',N'}(N'')=\delta'_{M',N'}(N'')
 =\delta_{M'',N''}(M')=\delta'_{M'',N''}(M')=0,\\
&\delta_{M',N'}(M'')=\delta'_{M',N'}(M'')
 =\delta_{M'',N''}(N')=\delta'_{M'',N''}(N')=0.
\end{aligned}
\end{equation}
By Theorem~\ref{codim1}, there are exact sequences
$$
0\to Z'\xrightarrow{f'}Z'\oplus M'\to N'\to 0\quad\text{and}\quad
0\to Z''\xrightarrow{f''}Z''\oplus M''\to N''\to 0
$$
in $\mod A$ such that the modules $Z'$ and $Z''$ are
indecomposable and
\begin{equation} \label{2outof3}
\delta'_{M',N'}(Z'\oplus M')=\delta'_{M'',N''}(Z''\oplus M'')=0.
\end{equation}
Observe that the homomorphisms $f'$ and $f''$ belong to
$\rad(\mod A)$, as they are not sections and $Z'$ and $Z''$ are
indecomposable modules.
Using \eqref{1outof3} and applying twice Lemma~\ref{degmix} we get
\begin{equation} \label{3outof3}
\delta'_{M',N'}(Z'')=\delta'_{M'',N''}(Z')=0.
\end{equation}
Let $M=M'\oplus M''$, $N=N'\oplus N''$ and $Z=Z'\oplus Z''$.
Taking a direct sum of the above exact sequences we obtain
an exact sequence of the form
$$
0\to Z\to Z\oplus M\to N\to 0.
$$
Applying \eqref{1outof3}, \eqref{2outof3} and \eqref{3outof3}
yields
$$
\delta'_{M,N}(Z\oplus M)=\delta'_{M',N'}(Z'\oplus M'\oplus Z''
 \oplus M'')+\delta'_{M'',N''}(Z'\oplus M'\oplus Z''\oplus M'')
=0.
$$
Hence $\Sing(M,N)=\Reg$, by Proposition~\ref{standcriterion}.
\end{proof}

We shall need the following result proved by Bongartz
in \cite[Theorem 5]{Bmin}.

\begin{prop} \label{fromUtoV}
Let $U$, $V$ and $M$ be modules such that
$\CO_{U\oplus V}\subseteq\ov{\CO}_M$ and
$\delta'_{M,U\oplus V}(U)=\delta_{M,U\oplus V}(V)=0$.
Then there is an exact sequence in $\mod A$ of the form
$0\to U\to M\to V\to 0$.
\end{prop}

\begin{prop} \label{ULV}
Let $M$ and $N$ be disjoint modules such that
$\CO_N\subseteq\ov{\CO}_M$.
Assume that $N\simeq U\oplus L\oplus V$ for some modules $U$, $L$
and $V$ such that
\begin{equation} \label{eqULV}
\begin{aligned}
\delta_{M,N}(U)&=1,&\delta_{M,N}(L)&=1,&\delta_{M,N}(V)&=0,
 &\delta_{M,N}(M)=0,\\
\delta'_{M,N}(U)&=0,&\delta'_{M,N}(L)&=1,&\delta'_{M,N}(V)&=1,
 &\delta'_{M,N}(M)=0.
\end{aligned}
\end{equation}
Then $\Sing(M,N)=\Reg$.
\end{prop}

\begin{proof}
Applying Theorem~\ref{specialseq} we get an exact sequence
$$
\sigma:\quad
0\to Z\xrightarrow{f}Z\oplus M\xrightarrow{g}N\to 0
$$
in $\mod A$ such that $f$ belongs to $\rad(\mod A)$.
Since $\delta'_{M,N}(U)=0$ and the modules $M$ and $U$ are
disjoint, then $Z\simeq U\oplus Y$ and there is an exact sequence
$$
\tau:\quad 0\to Z\xrightarrow{f'}Y\oplus M\to L\oplus V\to 0
$$
in $\mod A$ for some module $Y$ and some homomorphism $f'$ in
$\rad(\mod A)$, by Lemma~\ref{splitdirectsum}.
Taking a pushout of the sequence $\tau$ under a retraction
$\pi:Z\to U$ leads to the following commutative diagram with
exact rows and columns
$$
\xymatrix{
&0\ar[d]&0\ar[d]\\
&Y\ar@{=}[r]\ar[d]&Y\ar[d]\\
0\ar[r]&Z\ar[r]^-{f'}\ar[d]_\pi&Y\oplus M\ar[r]\ar[d]
 &L\oplus V\ar[r]\ar@{=}[d]&0\\
0\ar[r]&U\ar[r]\ar[d]&W\ar[r]\ar[d]&L\oplus V\ar[r]&0.\\
&0&0
}
$$
Applying Corollary~\ref{UMVdeg} and Theorem~\ref{specialseq}
to the exact sequences
$$
\varepsilon:\;
0\to U\xrightarrow{\alpha}W\xrightarrow{\beta}L\oplus V\to 0
\quad\text{and}\quad
0\to Y\to Y\oplus M\to W\to 0
$$
we get that $\CO_N\subseteq\ov{\CO}_W$ and
$\CO_W\subseteq\ov{\CO}_M$.
We conclude from \eqref{eqULV} the equality
$\delta'_{M,N}(U\oplus M)=0$.
Therefore if $W\simeq M$ then $\Sing(M,N)=\Reg$,
by Corollary~\ref{corcriterion} applied to the sequence
$\varepsilon$.
Thus we may assume that $W\not\simeq M$.
Since $f'$ belongs to $\rad(\mod A)$ then the retraction
$\pi$ does not factor through $f'$ and consequently the
exact sequence $\varepsilon$ does not split.
This implies that $W\not\simeq U\oplus L\oplus V\simeq N$,
by Lemma~\ref{whensplits}.
Therefore $\dim\CO_N<\dim\CO_W$ as well as $\dim\CO_W<\dim\CO_M$.
Since $\dim\CO_M-\dim\CO_N=2$ then
$$
\dim\CO_M-\dim\CO_W=1\quad\text{and}\quad
\dim\CO_W-\dim\CO_N=1.
$$
Applying Theorem~\ref{codim1} we get
\begin{equation} \label{help1}
\delta_{M,W}(W)=\delta'_{M,W}(W)=\delta_{W,N}(N)=\delta'_{W,N}(N)
=1,\quad\delta'_{W,N}(W)=0.
\end{equation}
Consequently
$$
1=\delta'_{W,N}(N)\geq\delta'_{W,N}(L\oplus V)=
\delta'_\varepsilon(L\oplus V)>0,
$$
by Lemma~\ref{whensplits}.
Thus
$$
\delta'_{W,N}(L)+\delta'_{W,N}(V)=1,
$$
which gives two possibilities.
\bigskip

\noindent\textbf{Case 1:}
$\delta'_{W,N}(L)=1$ and $\delta'_{W,N}(V)=0$.
\bigskip

Then $\delta'_\varepsilon(V)=0$, $W\simeq V\oplus W'$ and there
is an exact sequence
$$
\varepsilon':\quad 0\to U\to W'\to L\to 0
$$
in $\mod A$ for some module $W'$, by Lemma~\ref{splitdirectsum}.

It follows from \eqref{eqULV} and \eqref{help1} that
\begin{equation} \label{mideq}
\begin{aligned}
\delta'_{M,W}(W')&=\delta'_{M,W}(W)-\delta'_{M,W}(V)=1-
 (\delta'_{M,N}(V)-\delta'_{W,N}(V))=0,\\
\delta_{M,W}(V)&=\delta_{M,N}(V)-\delta_{W,N}(V)\leq
 \delta_{M,N}(V)=0.
\end{aligned}
\end{equation}
Hence $\delta_{M,W}(V)=0$ and there is an exact sequence
$$
\eta:\quad 0\to W'\to M\to V\to 0
$$
in $\mod A$, by Proposition~\ref{fromUtoV}.
It follows from \eqref{help1} that $\delta'_{W,N}(W')=0$.
Consequently, by \eqref{eqULV} and \eqref{mideq},
$$
\delta'_{M,N}(U\oplus W'\oplus M)=\delta'_{M,N}(W')
=\delta'_{M,W}(W')+\delta'_{W,N}(W')=0.
$$
Taking a direct sum of the sequences $\varepsilon'$, $\eta$
and $0\to 0\to U\xrightarrow{1_U}U\to 0$ gives an exact
sequence of the form
$$
0\to U\oplus W'\to U\oplus W'\oplus M\to N\to 0.
$$
Then $\Sing(M,N)=\Reg$, by Proposition~\ref{standcriterion}
applied for $Z=U\oplus W'$.
\bigskip

\noindent\textbf{Case 2:}
$\delta'_{W,N}(L)=0$ and $\delta'_{W,N}(V)=1$.
\bigskip

Then $\delta'_\varepsilon(L)=0$, $W\simeq L\oplus W''$ and there
is an exact sequence
$$
\varepsilon'':\quad 0\to U\to W''\to V\to 0
$$
in $\mod A$ for some module $W''$, by Lemma~\ref{splitdirectsum}.
In particular $U\oplus V\not\simeq W''$ and
$\CO_{U\oplus V}\subseteq\ov{\CO}_{W''}$,
by Corollary~\ref{UMVdeg}.
Applying Lemma~\ref{whensplits} to the sequence $\varepsilon$
yields $\delta_{W,N}(U)=\delta_\varepsilon(U)>0$.
Consequently
$$
\delta_{W,N}(L)=\delta_{W,N}(N)-\delta_{W,N}(U\oplus V)\leq
\delta_{W,N}(N)-\delta_{W,N}(U)\leq\delta_{W,N}(N)-1.
$$
It follows from \eqref{eqULV} and \eqref{help1} that
$\delta_{W,N}(N)-1=0$, $\delta_{W,N}(L)=0$ and
\begin{equation} \label{twozerodeltas}
\begin{aligned}
\delta_{M,W}(W'')&=\delta_{M,W}(W)-\delta_{M,W}(L)
 =1-(\delta_{M,N}(L)-\delta_{W,N}(L))=0,\\
\delta'_{M,W}(W'')&=\delta'_{M,W}(W)-\delta'_{M,W}(L)
 =1-(\delta'_{M,N}(L)-\delta'_{W,N}(L))=0.
\end{aligned}
\end{equation}
Let $Y$ be an indecomposable direct summand of $W''$.
Then $\delta_{M,W}(Y)=\delta'_{M,W}(Y)=0$ and
$\mu(M,Y)\geq\mu(W,Y)>0$, by Lemma~\ref{ostra}.
This implies that $M\simeq W''\oplus M'$ for some module $M'$
not isomorphic to $L$.
Furthermore $\CO_L\subset\ov{\CO}_{M'}$,
by \eqref{twozerodeltas} and Theorem~\ref{cancel}.
Applying Proposition~\ref{twodeg}, we get
$$
\Sing(M,N)=\Sing(W''\oplus M',(U\oplus V)\oplus L)=\Reg.
$$
This finishes the proof of Proposition~\ref{ULV}.
\end{proof}
\bigskip

\noindent
\textit{Proof of Theorem~\ref{aux2}.}
We decompose $N=N_1\oplus\cdots\oplus N_s$, where
$N_i$ is an indecomposable module for $i=1,\ldots,s=s(N)$.
Our assumptions and Lemma~\ref{formulacodim} imply that
$[N,N]-[M,M]=2$.
Therefore
\begin{equation} \label{deltywspolne}
\begin{aligned}
2&=\delta_{M,N}(M)+\sum_{i=1}^s\delta'_{M,N}(N_i)
 =\delta'_{M,N}(M)+\sum_{i=1}^s\delta_{M,N}(N_i),\\
4&=(\delta_{M,N}(M)+\delta'_{M,N}(M))+\sum_{i=1}^s\left(
 \delta_{M,N}(N_i)+\delta'_{M,N}(N_i)\right).
\end{aligned}
\end{equation}
Since the modules $M$ and $N$ are disjoint then
$\mu(M,N_i)=0$ and consequently
\begin{equation} \label{deltydlaNi}
\delta_{M,N}(N_i)+\delta'_{M,N}(N_i)\geq 1,
\qquad i=1,\ldots,s,
\end{equation}
by Lemma~\ref{ostra}.
This implies that $s\leq 4$.
Recall that $s\geq 3$, by our assumptions.
Hence
\begin{equation} \label{deltydlaM}
\delta_{M,N}(M)+\delta'_{M,N}(M)\leq 1.
\end{equation}

Let $U$ and $V$ be the direct sums of the modules $N_i$ such that
$\delta'_{M,N}(N_i)=0$ and $\delta_{M,N}(N_i)=0$, respectively.
Then $\delta'_{M,N}(U)=0$ and $\delta_{M,N}(V)=0$.
It follows from \eqref{deltywspolne} and \eqref{deltydlaNi} that
$N\simeq U\oplus V\oplus L$, where either $L=0$, or $L=N_j$ for
some $j\leq s$ and the equalities \eqref{eqULV} hold.
We get $\Sing(M,N)=\Reg$ in the latter case, by
Proposition~\ref{ULV}.
Therefore we may assume that $L=0$, or equivalently,
$N\simeq U\oplus V$.
Then there is an exact sequence
$$
0\to U\to M\to V\to 0
$$
in $\mod A$, by Proposition~\ref{fromUtoV}.
Furthermore, \eqref{deltydlaM} implies that
$\delta_{M,N}(M)=0$ or $\delta'_{M,N}(M)=0$.
Hence $\Sing(M,N)=\Reg$, by Corollary~\ref{corcriterion}.
This finishes the proof of Theorem~\ref{aux2}.
\qed

\section{Path algebras of Dynkin quivers}
\label{proofmain}

Throughout the section, $A$ is the path algebra of a Dynkin
quiver.
We shall need some special properties of modules over such
algebra $A$ described in the following three lemmas, in order
to prove Theorem~\ref{main}.
The first lemma follows from \cite{Gdyn} and the second one
follows from \cite[Lemma 5]{Bmin}.

\begin{lem} \label{fundDynkin}
There are only finitely many isomorphism classes
of indecomposable modules.
Moreover, for each indecomposable module $Y$,
$$
\End_A(Y)=\left\{t\cdot 1_Y;\;t\in k\right\}.
$$
\end{lem}

\begin{lem} \label{notwice}
Let $M$ and $N$ be disjoint modules such that
$\CO_N\subset\ov{\CO}_M$ and $\dim\CO_M-\dim\CO_N=1$.
Then the inequality $\mu(M,Y)\leq 1$ holds for any
indecomposable module $Y$.
\end{lem}

\begin{lem} \label{specialUV}
Let $M$ and $N$ be disjoint modules with $\CO_N\subset\ov{\CO}_M$.
Then there are indecomposable direct summands $U$ and $V$ of $N$
such that
$$
\delta_{M,N}(U)>0,\;\delta'_{M,N}(U)=0
\quad\text{and}\quad
\delta_{M,N}(V)=0,\;\delta'_{M,N}(V)>0.
$$
\end{lem}

\begin{proof}
A complete set $\ind A$ of pairwise nonisomorphic indecomposable
modules is finite, by Lemma~\ref{fundDynkin}.
Moreover there is a partial order $\preceq$ on $\ind A$ such
that $[X,Y]>0$ implies $X\preceq Y$ for any modules $X$ and $Y$
in $\ind A$.
Applying Theorem~\ref{specialseq} we get an exact sequence
$$
\eta:\quad 0\to N\to M\oplus Z'\to Z'\to 0
$$
in $\mod A$.
Then $\delta_{M,N}(N)=\delta_\eta(N)>0$, by
Lemma~\ref{whensplits}.
Hence there is a $\preceq$-minimal $U\in\ind A$ with the property
$\delta_{M,N}(U)>0$.
Then $\mu(N,U)>0$, by \cite[Lemma 3.1]{Bext}.
Moreover, using the Auslander-Reiten formula mentioned in the
proof of \cite[Lemma 3.1]{Bext}, we get that $\delta'_{M,N}(U)=0$.
Dually we get an appropriate module $V$.
\end{proof}

\begin{prop} \label{longprop}
Let $\sigma:\;0\to U\xrightarrow{f}M\xrightarrow{g}V\to 0$ be
an exact sequence in $\mod A$ such that the modules $M$ and
$N=U\oplus V$ are disjoint and
\begin{equation} \label{assumpt}
\begin{aligned}
\delta_\sigma(U)&=1,&\delta_\sigma(M)&=1,&\delta_\sigma(V)&=0,\\
\delta'_\sigma(U)&=0,&\delta'_\sigma(M)&=1,&\delta'_\sigma(V)&=1.
\end{aligned}
\end{equation}
Then $\Sing(M,N)=\Reg$.
\end{prop}

\noindent
\textit{Proof of Proposition~\ref{longprop}.}
The equality $\delta_\sigma(M)=1$ implies that
$M=M_1\oplus M'$ for an indecomposable module $M_1$
and a module $M'$ such that
\begin{equation} \label{eqM1}
\delta_\sigma(M_1)=1\quad\text{and}\quad\delta_\sigma(M')=0.
\end{equation}
We divide the proof into several steps.

\begin{step} \label{step1}
There are nonsplittable exact sequences in $\mod A$ of the form
$$
\sigma_1:\;0\to U\xrightarrow{\psmatrix{f\\ h}}M\oplus M_1
 \xrightarrow{(h',-f')}X\to 0,\quad
\sigma_2:\;0\to M_1\xrightarrow{f'}X\xrightarrow{g'}V\to 0.
$$
\end{step}

\begin{proof}
Since $\delta_\sigma(M_1)>0$ then there is a homomorphism
$h:U\to M_1$ which does not factor through $f$, by
Lemma~\ref{seqineq}.
Taking a pushout of $\sigma$ under $h$ leads to the following
commutative diagram with exact rows
$$
\xymatrix{
0\ar[r]&U\ar[r]^f\ar[d]_h&M\ar[r]^g\ar[d]^{h'}&V\ar[r]\ar@{=}[d]
 &0\\
0\ar[r]&M_1\ar[r]^{f'}&X\ar[r]^{g'}&V\ar[r]&0,
}
$$
This gives the exact sequences $\sigma_1$ and $\sigma_2$.
The sequence $\sigma_2$ does not split, by our construction.
Since the modules $U$ and $M\oplus M_1$ are disjoint,
the sequence $\sigma_1$ does not split as well.
\end{proof}

\begin{step} \label{step2}
The following equalities hold:
\begin{equation} \label{eqforstep2}
\begin{aligned}
\delta_{\sigma_1}(U)&=1,&\delta_{\sigma_1}(M)&=0,
 &\delta'_{\sigma_1}(U)&=0,\\
\delta_{\sigma_2}(V)&=0,&\delta'_{\sigma_2}(U)&=0,
 &\delta'_{\sigma_2}(V)&=1.
\end{aligned}
\end{equation}
\end{step}

\begin{proof}
Since the sequences $\sigma_1$ and $\sigma_2$ do not split
then the integers $\delta_{\sigma_1}(U)$, $\delta_{\sigma_2}(M_1)$
and $\delta'_{\sigma_2}(V)$ are positive,
by Lemma~\ref{whensplits}.
Hence the claim follows from \eqref{assumpt}, \eqref{eqM1}
and the equalities
$$
\delta_\sigma(Y)=\delta_{\sigma_1}(Y)+\delta_{\sigma_2}(Y)
\quad\text{and}\quad
\delta'_\sigma(Y)=\delta'_{\sigma_1}(Y)+\delta'_{\sigma_2}(Y),
$$
for any module $Y$.
\end{proof}

\begin{step} \label{step3}
$\delta_{\sigma_1}(X)=0$.
\end{step}

\begin{proof}
Let $\tilde{M}=M\oplus M_1$.
The sequence $\sigma_1$ induces the following commutative diagram
with exact rows and columns
$$
\xymatrix{
&0\ar[d]&0\ar[d]&0\ar[d]\\
0\ar[r]&\Hom_A(X,U)\ar[r]\ar[d]&\Hom_A(X,\tilde{M})\ar[r]\ar[d]
 &\Hom_A(X,X)\ar[d]\\
0\ar[r]&\Hom_A(\tilde{M},U)\ar[r]\ar[d]
 &\Hom_A(\tilde{M},\tilde{M})\ar[r]\ar[d]^\alpha
 &\Hom_A(\tilde{M},X)\ar[d]^\gamma\\
0\ar[r]&\Hom_A(U,U)\ar[r]&\Hom_A(U,\tilde{M})\ar[r]^\beta
 &\Hom_A(U,X).
}
$$
Since $\delta_{\sigma_1}(\tilde{M})=\delta'_{\sigma_1}(U)=0$,
then the homomorphisms $\alpha$ and $\beta$ are surjective.
Hence $\gamma$ is also surjective, which implies that
$\delta_{\sigma_1}(X)=0$.
\end{proof}

\begin{step} \label{step4}
$\delta'_{\sigma_2}(M)=0$.
\end{step}

\begin{proof}
Suppose that $\delta'_{\sigma_2}(M)\geq 1$.
Since
$1=\delta_\sigma(M)=\delta_{\sigma_1}(M)+\delta_{\sigma_2}(M)$,
then
$$
\delta'_{\sigma_1}(M)=0\quad\text{and}\quad
\delta'_{\sigma_1}(M_1)=0,
$$
as $M_1$ is a direct summand of $M$.
Observe that
$$
\delta_{\sigma_1}(U)-\delta_{\sigma_1}(M\oplus M_1)
+\delta_{\sigma_1}(X)=
\delta'_{\sigma_1}(U)-\delta'_{\sigma_1}(M\oplus M_1)
+\delta'_{\sigma_1}(X).
$$
Applying \eqref{eqforstep2} and Step~\ref{step3} we get that
$\delta'_{\sigma_1}(X)=1$.
Then $X=X_1\oplus X'$ for an indecomposable module $X_1$ and
a module $X'$ such that
$$
\delta'_{\sigma_1}(X_1)=1\quad\text{and}\quad
\delta'_{\sigma_1}(X')=0.
$$
Let $\varphi:X'\to X$ be a section.
Hence $\varphi=h'\tilde{h}-f'\tilde{f}$ for some homomorphisms
$\tilde{h}:X'\to M$ and $\tilde{f}:X'\to M_1$, by
Lemma~\ref{seqineq} applied to the sequence $\sigma_1$.
Since the sequence $\sigma_2$ does not split and the module $M_1$
is indecomposable, then $f'$ belongs to $\rad(\mod A)$.
Thus $f'\tilde{f}$ belongs to $\rad(\mod A)$ and $h'\tilde{h}$ is
a section.
Consequently $\tilde{h}$ is also a section.
Applying Lemma~\ref{howtosplitiso} to $\sigma_1$ we get that
$M\simeq X'\oplus M''$ and there is an exact sequence
$$
\tau:\quad 0\to U\to M''\oplus M_1\to X_1\to 0
$$
in $\mod A$ for some module $M''$.
The modules $U$ and $M''\oplus M_1$ are disjoint, by our
assumptions.
The modules $X_1$ and $M''\oplus M_1$ are also disjoint,
since $X_1$ is indecomposable, $\delta'_{\sigma_1}(X_1)>0$
and $\delta'_{\sigma_1}(M''\oplus M_1)=0$.
Observe that
$$
\dim\CO_{M''\oplus M_1}-\dim\CO_{U\oplus X_1}
=\delta_{\sigma_1}(U\oplus X_1)
+\delta'_{\sigma_1}(M''\oplus M_1)=1,
$$
by \eqref{eqforstep2} and Step~\ref{step2}.
Hence $\mu(M'',M_1)=0$, by Lemma~\ref{notwice}.
Since $M_1\oplus M'$ is isomorphic to $X'\oplus M''$ then
$\mu(X,M_1)\geq\mu(X',M_1)\geq 1$ and $X\simeq M_1\oplus X''$
for some module $X''$.
Hence, up to an isomorphism, the sequence $\sigma_2$ has the form
$$
0\to M_1\xrightarrow{f'=\psmatrix{\alpha_1\\ \alpha_2}}
M_1\oplus X''\xrightarrow{g'=(\beta_1,\beta_2)}V\to 0.
$$
Since the endomorphism $\alpha_1\in\End_A(M_1)$ belongs to
$\rad(\mod A)$ and $M_1$ is an indecomposable module, then
$\alpha_1=0$, by Lemma~\ref{fundDynkin}.
Observe that
$$
\Ker(\beta_1)\subseteq\Ker(g')\cap M_1\quad\text{and}\quad
\Ker(g')=\im(f')\subseteq X''.
$$
Therefore the homomorphism $\beta_1$ is injective and
$\im(\beta_1)\cap\im(\beta_2)=\{0\}$.
Thus $\im(\beta_1)$ is a direct summand of $V$, as $g'$ is
surjective.
Consequently the homomorphism $\beta_1:M_1\to V$ is a section,
which is impossible as $M_1$ and $V$ are disjoint modules.
\end{proof}

\begin{step} \label{step5}
$\delta_\sigma(X)=0$.
\end{step}

\begin{proof}
Observe that
$$
\delta_\sigma(M_1)-\delta_\sigma(X)+\delta_\sigma(V)
=\delta'_{\sigma_2}(U)-\delta'_{\sigma_2}(M)
+\delta'_{\sigma_2}(V).
$$
Hence the claim follows from \eqref{assumpt}, \eqref{eqM1},
\eqref{eqforstep2} and Step~\ref{step4}.
\end{proof}

\begin{step} \label{step6}
There is an exact sequence
$\sigma_3:\;0\to U\to X\oplus M'\to V\oplus V\to 0$.
\end{step}

\begin{proof}
Since $M=M_1\oplus M'$ then the sequence $\sigma$ has the form
$$
0\to U\xrightarrow{\psmatrix{f_1\\ f_2}}M_1\oplus M'
\xrightarrow{(g_1,g_2)}V\to 0.
$$
We get from \eqref{eqforstep2} the equality
$\delta_{\sigma_2}(V)=0$.
Hence any homomorphism from $M_1$ to $V$ factors through $f'$,
by Lemma~\ref{seqineq}.
Thus $g_1=jf'$ for some homomorphism $j:X\to V$.
It is easy to check that the sequence
$$
0\to U\xrightarrow{\psmatrix{f'f_1\\ f_2}}X\oplus M'
\xrightarrow{\psmatrix{g'&0\\ j&g_2}}V\oplus V\to 0
$$
is exact.
\end{proof}

We shall consider the $k$-functor $\CE_{M,N}(-,-)$ defined in
Section~\ref{smooth}.

\begin{step} \label{step7}
$\dim_k\CE_{M,N}(V,U)\leq 2$.
\end{step}

\begin{proof}
We know that
$\delta_{M,N}(X\oplus M')=\delta_\sigma(X)+\delta_\sigma(M')=0$,
by \eqref{eqM1} and Step~\ref{step5}.
Applying Corollary~\ref{Efunctor} we get that $\CE_{M,N}(V,U)$ is
contained in the kernel of the last map in the following long
exact sequence induced by $\sigma_3$:
\begin{multline*}
0\to\Hom_A(V,U)\to\Hom_A(V,X\oplus M')\to\Hom_A(V,V\oplus V)\to\\
\to\Ext^1_A(V,U)\to\Ext^1_A(V,X\oplus M').
\end{multline*}
Consequently
$$
\dim_k\CE_{M,N}(V,U)\leq\delta'_{\sigma_3}(V)
=\delta'_\sigma(V)+\delta'_{\sigma_2}(V)=1+1=2,
$$
by \eqref{assumpt} and \eqref{eqforstep2}.
\end{proof}

\begin{step} \label{step8}
$\dim_k\CE_{M,N}(N,N)\leq[N,N]-[M,M]$.
\end{step}

\begin{proof}
Let $Y$ be a module.
We know that $\delta_{M,N}(V)=\delta_\sigma(V)=0$ and
$\delta'_{M,N}(U)=\delta'_\sigma(U)=0$, by \eqref{assumpt}.
Then $\CE_{M,N}(Y,V)$ is contained in the kernel of
$\Ext^1_A(Y,1_V)$ and $\CE_{M,N}(U,Y)$ is contained in the
kernel of $\Ext^1_A(1_U,Y)$, by Corollary~\ref{Efunctor}.
Hence $\CE_{M,N}(Y,V)=0$ and $\CE_{M,N}(U,Y)=0$.
Consequently
\begin{multline*}
\CE_{M,N}(N,N)\simeq\CE_{M,N}(U\oplus V,U\oplus V)\\
\simeq\CE_{M,N}(U,U)\oplus\CE_{M,N}(U,V)
\oplus\CE_{M,N}(V,U)\oplus\CE(V,V)\simeq\CE_{M,N}(V,U).
\end{multline*}
Therefore the claim follows from Step~\ref{step7} and the
equalities
$$
[N,N]-[M,M]=\delta_\sigma(U)+\delta_\sigma(V)+\delta'_\sigma(M)
=1+0+1=2.
$$
\end{proof}

Step~\ref{step8} together with Proposition~\ref{gencriterion}
imply that $\Sing(M,N)=\Reg$, which finishes the proof
of Proposition~\ref{longprop}.
\qed
\bigskip

\noindent
\textit{Proof of Theorem~\ref{main}.}
Let $M$ be a module.
It follows from Lemma~\ref{fundDynkin} that $\ov{\CO}_M$
contains only finitely many orbits.
Thus it suffices to show that $\Sing(M,N)=\Reg$ for any module
$N$ such that $\CO_N\subset\ov{\CO}_M$ and
$$
c:=\dim\CO_M-\dim\CO_N\in\{1,2\}.
$$
If $c=1$, then the claim follows from Theorem~\ref{codim1}.
Therefore we may assume that $c=2$.
Applying Theorem~\ref{aux1} we reduce the problem to the case
when the modules $M$ and $N$ are disjoint.
Then $N\simeq U\oplus V\oplus L$, where $U$ and $V$ are
indecomposable modules such that
$$
\delta_{M,N}(U)>0,\;\delta'_{M,N}(U)=0
\quad\text{and}\quad
\delta_{M,N}(V)=0,\;\delta'_{M,N}(V)>0,
$$
by Lemma~\ref{specialUV}.
Applying Theorem~\ref{aux2} we may assume that $L=0$ and
$N\simeq U\oplus V$.
Hence there is an exact sequence
$$
\sigma: 0\to U\to M\to V\to 0,
$$
by Proposition~\ref{fromUtoV}.
If $\delta_\sigma(M)=0$ or $\delta'_\sigma(M)=0$ then
$\Sing(M,N)=\Reg$, by Corollary~\ref{corcriterion}.
Therefore we may assume that the integers $\delta_\sigma(M)$
and $\delta'_\sigma(M)$ are positive.
On the other hand, by Lemma~\ref{formulacodim},
\begin{align*}
2=[N,N]-[M,M]
 &=\delta_\sigma(U)+\delta_\sigma(V)+\delta'_\sigma(M)\\
 &=\delta'_\sigma(U)+\delta'_\sigma(V)+\delta_\sigma(M),
\end{align*}
which implies that the equalities \eqref{assumpt} hold.
Thus $\Sing(M,N)=\Reg$, by Proposition~\ref{longprop}.
This finishes the proof of Theorem~\ref{main}.
\qed

\bigskip

\noindent
Grzegorz Zwara\\
Faculty of Mathematics and Computer Science\\
Nicolaus Copernicus University\\
Chopina 12/18\\
87-100 Toru\'n\\
Poland\\
E-mail: gzwara@mat.uni.torun.pl
\end{document}